\documentclass{amsart}

\usepackage{indentfirst}
\usepackage{amssymb}
\usepackage{amsmath}
\usepackage{graphics}
\usepackage{epsfig}
\usepackage{xcolor}
\usepackage{MnSymbol,wasysym}
\usepackage{mathtools}

\usepackage[utf8]{inputenc}
\usepackage{latexsym}
\usepackage{graphicx}

\usepackage[normalem]{ulem}
\usepackage{soul}

\usepackage{scalerel,stackengine}
\stackMath
\newcommand\reallywidecheck[1]{%
\savestack{\tmpbox}{\stretchto{%
  \scaleto{%
    \scalerel*[\widthof{\ensuremath{#1}}]{\kern-.6pt\bigwedge\kern-.6pt}%
    {\rule[-\textheight/2]{1ex}{\textheight}}
  }{\textheight}%
}{0.5ex}}%
\stackon[1pt]{#1}{\scalebox{-1}{\tmpbox}}%
}

\catcode`@=11 \@addtoreset{equation}{section}
\renewcommand\theequation{\thesection.\@arabic\c@equation}
\catcode`@=12

\newcommand{\RR}{\mathbb{R}}

\expandafter\chardef\csname pre amssym.def
at\endcsname=\the\catcode`\@ \catcode`\@=11
\def\undefine#1{\let#1\undefined}
\def\newsymbol#1#2#3#4#5{\let\next@\relax
 \ifnum#2=\@ne\let\next@\msafam@\else
 \ifnum#2=\tw@\let\next@\msbfam@\fi\fi
 \mathchardef#1="#3\next@#4#5}
\def\mathhexbox@#1#2#3{\relax
 \ifmmode\mathpalette{}{\m@th\mathchar"#1#2#3}%
 \else\leavevmode\hbox{$\m@th\mathchar"#1#2#3$}\fi}
\def\hexnumber@#1{\ifcase#1 0\or 1\or 2\or 3\or 4\or 5\or 6\or 7\or 8\or
 9\or A\or B\or C\or D\or E\or F\fi}

\font\teneufm=eufm10 \font\seveneufm=eufm7 \font\fiveeufm=eufm5
\newfam\eufmfam
\textfont\eufmfam=\teneufm \scriptfont\eufmfam=\seveneufm
\scriptscriptfont\eufmfam=\fiveeufm

\catcode`\@=\csname pre amssym.def at\endcsname

\newcommand{\eqn}{\begin{eqnarray}}
\newcommand{\een}{\end{eqnarray}}

\newtheorem {Theorem}  {Theorem}

\numberwithin{Theorem}{section}
\newtheorem{Lemma}[Theorem]{Lemma}

\theoremstyle{Remark}
\newtheorem{Remark}[Theorem]{Remark}

\newcommand{\NN}{{\mathbb N}}
\newcommand{\ZZ}{{\mathbb Z}}

\newcommand{\RRn}{{\mathbb R}^n}

\renewcommand{\a}{\alpha}

\renewcommand{\t}{\tau}

\newcommand{\ff}{\mbox{\boldmath $f$}}
 
\newcommand{\uu}{\mbox{\boldmath $u$}}

\newcommand{\BB}{\mbox{\boldmath $B$}}
\newcommand{\bb}{\mbox{\boldmath $b$}}

\begin{document}

\title[$\dot{H}^{m}$ estimates for parabolic equations]{Upper and lower $\dot{H}^{m}$ estimates
for solutions to parabolic equations}

\author[R.H. Guterres]{Robert H. Guterres}
\address[R.H. Guterres]{Departamento de Matemática. Universidade Federal de Pernambuco, PE 50740-560, Recife. Brazil.}
\email{rguterres.mat@gmail.com}

\author[C. J. Niche]{C\'esar J. Niche}
\address[C.J. Niche]{Departamento de Matem\'atica Aplicada, Instituto de Matem\'atica. Universidade Federal do Rio de Janeiro, CEP 21941-909, Rio de Janeiro - RJ, Brazil.}
\email{cniche@im.ufrj.br}

\author[C.F. Perusato]{Cilon F. Perusato}
\address[C.F. Perusato]{Departamento de Matem\'atica. Universidade Federal de Pernambuco, PE 50740-560, Recife. Brazil.}
\email{cilon.perusato@ufpe.br}

\author[P.R. Zingano]{Paulo R. Zingano}
\address[P.R. Zingano]{Departamento de Matemática Pura e Aplicada. Universidade Federal do Rio Grande do Sul, RS 91509-900, Porto Alegre. Brazil.}
\email{paulo.zingano@ufrgs.br }

\thanks{R. H. Guterres acknowledges support  form FACEPE Grant BFP-0067-1.01/19. C.J. Niche acknowledges support from Bolsa PQ CNPq - 308279/2018-2  and PROEX - CAPES. C.J. Niche and C.F. Perusato acknowledge support from PRONEX-FAPERJ ``Equa\c{c}\~oes Diferenciais Parciais N\~ao Lineares e Aplica\c{c}\~oes''. C. F. Perusato was partially supported by CAPES - PRINT - 88881.311964/2018-01. He is also grateful for the warm hospitality during his visit at the Universidade Federal do Rio de Janeiro, where this
work was started.  } 



\date{\today}

\begin{abstract}
In this article we prove results concerning upper and lower decay estimates for homogeneous Sobolev norms of solutions to a rather general family of parabolic equations. Following the ideas of Kreiss, Hagstrom, Lorenz and Zingano, we use eventual regularity of solutions to directly work with smooth solutions in physical space, bootstrapping decay estimates from the $L^2$ norm to higher order derivatives. Besides obtaining upper and lower bounds through this method, we also obtain reverse results: from higher order derivatives decay estimates, we deduce bounds for the $L^2$ norm. We use these general results to prove new decay estimates for some equations and to recover some well known results.    
\end{abstract}

\maketitle

\section{Introduction}

\subsection{Settings} The class of parabolic equations of the form 

\begin{equation}
\label{eqn:main-equation}
\uu_t + \mathbb{G}(\uu) = \nu \Delta \uu +  \ff, \qquad x \in \RR^n, t > 0,
\end{equation}
where $\uu$ is a vector field, $\mathbb{G}(\uu)$ is a nonlinear function of $\uu$, $\ff = \ff (x,t)$ is an external force and $\nu > 0$ is a constant, includes a large number of important equations on which many results concerning upper and lower estimates for evolution of homogeneous Sobolev norms $\dot{H} ^m (\RR ^n), m \in \NN$,  have been proved. As examples of equations that fit in the framework of (\ref{eqn:main-equation}) we mention the Navier-Stokes equations, advection-diffusion equations, MHD and MHD-like equations, micropolar and magneto-micropolar equations, amongst many others. Our main goal in this article is to prove general results providing upper and lower bounds for the decay rate of the $\dot{H} ^m (\RR ^n)$ norm, provided solutions $\uu$ have certain natural properties which are known to hold for  equations as those listed above. 

We describe now the setting in which we work. We assume first that (\ref{eqn:main-equation}) has  global weak solutions $ \uu(\cdot,t) \in C_{\tt w} \left( [0, \infty), L^{2}(\RR ^{n}) \right)$  which become eventually smooth and that $\ff$ also becomes eventually smooth,  i.e.  for some $t_{\ast} > 0$ 

\begin{align}
\label{eqn:evreg1} \ff, \uu  & \in  C^{\infty}(\RR ^{n} \times (t_{\ast}, \infty)), \\
\label{eqn:evreg2} \ff (\cdot,t), \uu(\cdot,t) & \in  C^{0}((t_{\ast}, \infty), H^{m}(\mathbb{R}^{n}) ^n),
\quad \, \forall \, m \geq 0.
\end{align}
In order to prove our results, we consider the following set of hypotheses:

\begin{itemize}
\item there exist constants $T_0, C_0, \alpha > 0$ such that for all $t > T_0$, 

\begin{equation} \tag{H1} \label{eqn:h1}
\Vert \uu(\cdot,t) \Vert _{L^2} \leq C_0  t^{- \alpha};
\end{equation}
\item for some $ \widehat{m} \geq 1 $ we have, for each $ 0 \leq m \leq \widehat{m} $, and some $ \tau_m > t_{\ast} $, 

\begin{equation} \tag{H2} \label{eqn:h2}
\big| \sum_{\ell_{1},...,\ell_{m}} \int_{\RR ^n} \langle D_{\ell_1} \cdots  D_{\ell_m} \uu (x,t),   D_{\ell_1} \cdots  D_{\ell_m} \mathbb{G}(\uu(x,t)) \rangle \,dx\, \big|   \leq g_{m} (t) \| D^{m+1} \uu( \cdot,t) \|^{2} _{L^2},
\end{equation}
for all $ t>\tau_m $, where the sum is over all indices $ 1 \leq \ell_{1}, \cdots, \ell_{m} \leq n $ and $g_m(t)$ is a continuous function that tends to zero when $t$ goes to infinity;
\item for some constant $ \beta>0 $ and each $ 0 \leq m \leq \widehat{m} $, there exist $ F_m, \sigma_m > 0$ such that 
\begin{equation}
\tag{H3} \label{eqn:h3}
\| D^{m} \ff (\cdot,t) \| _{L^2} \leq  F_m  t^{- \beta -  \frac{m}{2} },
\end{equation}
for all $ t>\sigma_m $;
\item for some constants $D_0, t_0, \eta$ and $ t > t_0$

\begin{equation}
\tag{H4} \label{eqn:h4}
\| \uu(\cdot,t) \| _{L^2}  \geq  D_0 t^{- \eta},
\end{equation}
where $\eta \geq \alpha$, for $\alpha$ as in (\ref{eqn:h1}).
\end{itemize}
We discuss now the assumption and hypotheses needed to prove our results. We first assume eventual regularity, i.e. (\ref{eqn:evreg1}) and (\ref{eqn:evreg2}). This is key to our method, as once the solution is smooth enough we directly work with equation (\ref{eqn:main-equation}) taking derivatives and obtaining energy equalities and inequalities, where we use the hypotheses to bootstrap estimates as in (\ref{eqn:h1}) and  (\ref{eqn:h4}) from lower to higher order derivatives through (\ref{eqn:h2}) (and vice versa) and rather simple calculations. Note that we do not assume a specific functional form for the nonlinearity $ \mathbb{G}(\uu(\cdot,\cdot))$, but we only ask for decay in time estimates on its derivatives through (\ref{eqn:h2}). For many equations of the form (\ref{eqn:main-equation}) these estimates can be naturally proved by induction with (\ref{eqn:h1}) and  (\ref{eqn:h4}) as base step. It is important to note that these estimates and computations all take place in physical space, see Remak \ref{comparison} below. When there is an external force acting on the equation, we need the independent estimate (\ref{eqn:h3}). We closely follow the ideas in Hagstrom, Lorenz, J. P. Zingano and P. Zingano \cite{MR4021907} and Kreiss, Hagstrom, Lorenz and P. Zingano \cite{MR1994780}.

\subsection{Results} We start now stating the main results in this article.

\begin{Theorem}
\label{upper-bound-decay}
Assume (\ref{eqn:h1}),  (\ref{eqn:h2}) and  (\ref{eqn:h3}) as above, with $\beta \geq \alpha + 1 $. Then, for every $ 1 \leq m \leq \widehat{m}$, there exist constants $ C_m > 0, \, T_m > t_{\ast} $ such that 

\begin{displaymath}
\|D^{m} \uu (\cdot,t) \| _{L^2} \leq \, C_m  \,\nu^{-\frac{m}{2}} \, t^{- \alpha -  \frac{m}{2} }, \qquad t > T_m.
\end{displaymath}
\mbox{ }\\
If $\beta > \alpha +1$, then $C_m = C_m \left(C_0, \alpha \right)$, and if $\beta = \alpha +1$, then  $C_m = C_m (C_0, \nu, \alpha,  F_{0}, $ $ \ldots, F_{m-1} )$. The time $T_m$ is such that  $T_m = T_m\,(C_0, \nu, \alpha, \beta,  F_{0}, \ldots,  F_{m},\tau_{0},\ldots, \tau_{m}, \sigma_{0},$ $ \ldots, \sigma_{m}, g_{0}, \ldots, g_{m}  ) $.
\end{Theorem}
Roughly speaking, every derivative contributes minus one half  to the decay rate of the solution $u$, which leads to the correct estimate provided the external force has fast enough decay. This result is well known for the Navier-Stokes equations, see Bae and Biswas \cite{MR3457535}, Benameur and Selmi \cite{MR2975158}, Biswas \cite{MR2964642}, Deng and Shang \cite{MR4220767}, Hagstrom, Lorenz, J.P. Zingano and P. Zingano \cite{MR4021907}, Jiu and Yu \cite{MR3385245}, Niche and M.E. Schonbek \cite{MR3355116}, Oliver and Titi \cite{MR1749867}, M.E. Schonbek \cite{MR1312701}, M.E. Schonbek and Wiegner \cite{MR1396285}, and references therein; for the MHD equations and related models, see  Ahn \cite{MR4103923}, Chae and M.E. Schonbek \cite{MR3097244}, de Souza, Melo and Zingano \cite{MR3898701}, Duan, Fukumoto and Zhao \cite{MR3985373}, Li , Zhu and Zhao \cite{MR3529623}, Melo, Perusato, Guterres and Nunes \cite{MR4128674}, Weng \cite{MR3460222}, \cite{MR3460237}, Wu, Yu and Tang \cite{MR3596669}, Zhao \cite{MR3850026}, \cite{MR4122514}, \cite{MR4102835}, Zhao and Zhu \cite{MR3823967}, \cite{MR3760328}, \cite{MR4054821} ; for families of dissipative equations, see Bae and Biswas \cite{MR3457535}, Biswas \cite{MR2964642}, Braz e Silva, Guterres, Perusato and P. Zingano \cite{superdesigualdade_submetido}, Niche and M.E. Schonbek \cite{MR3355116};  and for many other equations. 

\begin{Remark} Due to the rather general form of (\ref{eqn:main-equation}) and the natural hypotheses  (\ref{eqn:h1}),  (\ref{eqn:h2}) and  (\ref{eqn:h3}), we expect Theorem \ref{upper-bound-decay} to be widely applicable to other  equations. $\Box$
\end{Remark}

\begin{Remark}
Theorem \ref{upper-bound-decay} is also valid  in the case $ \alpha = 0 $, as it will become clear from  its derivation. Further generalizations are possible;  for example, if instead of (\ref{eqn:h1}) we assume that $\uu(t) = o(t^{-\alpha})$  as $ t \rightarrow \infty $, for some $ \alpha \geq 0 $, then we have that $ \| D^{m} \uu(\cdot,t) \| _{L^2} = o(t^{- \alpha -  \frac{m}{2} })$, for every $1 \leq m \leq \widehat{m} $. $\Box$
\end{Remark}

{\begin{Remark}
For a detailed and thorough survey of results concerning decay of solutions to the Navier-Stokes equations,  see Brandolese and M.E. Schonbek \cite{MR3916783}.
\end{Remark}

We  state now results concerning lower bounds for the decay of solutions.

\begin{Theorem}[The case $\eta = \alpha$]
\label{lower-bound-decay-first}
Assume (\ref{eqn:h1}), (\ref{eqn:h2}), (\ref{eqn:h3}) and (\ref{eqn:h4}) with $\eta = \alpha$ and $ \beta > \alpha + 1 $. Then, for every $ 1  \leq m \leq \widehat{m} $ we have that 

\begin{displaymath}
\| D^{m} \uu(\cdot,t) \| _{L^2} \geq\, D_m \, \nu^{-\frac{m}{2}}\, t^{- \alpha \,-\, \frac{m}{2} }, \quad \forall  \, t > t_m.
\end{displaymath}
The constants $D_m$ are such that $D_m = D_m (\alpha,  C_0, D_0)$, while $t_m > t_{\ast}$ is such that $t_m = t_m \,(C_0, D_0, \nu, \alpha, \beta,$ $t_0,T_0, F_{0}, \ldots, F_{m},\tau_{0},\ldots, \tau_{m},\sigma_{0}, \ldots, \sigma_{m}, g_{0}, \ldots, g_{m}  ) $. 
\end{Theorem}
\newpage
\begin{Theorem}[The case $\eta > \alpha$]
\label{lower-bound-decay-second}
Assume (\ref{eqn:h1}), (\ref{eqn:h2}), (\ref{eqn:h3}) and (\ref{eqn:h4}) with $\eta > \alpha$ and

\begin{displaymath}
\beta > 2 \eta \,-\, \alpha + \left( \frac{\eta}{\alpha} - 1 \right) \widehat{m} + 1.
\end{displaymath}
Then, setting $q = \frac{\eta}{\alpha}$,  we have that for every $ 1  \leq m \leq \widehat{m} $ 

\begin{displaymath}
\| D^{m} \uu(\cdot,t) \| _{L^2} \geq {D}_m\,\nu^{-\frac{m}{2}}\,  t^{-\eta - \frac{m q}{2}}, \quad \forall  t > t_m.
\end{displaymath}
Both ${D}_m$ and $t_m$ depend on the same parameters as in Theorem \ref{lower-bound-decay-first}, except that $ t_m $ depends also on $ \eta $.
\end{Theorem}

Theorem \ref{lower-bound-decay-second}  is a slight generalization of Theorem \ref{lower-bound-decay-first} and applies to the case where the lower bound for decay rate is faster than the upper bound (see Remark \ref{remark-decay-character}). The proofs of these results are based on the use of energy estimates for  $ \| D^{m-1} \uu(\cdot,t)\|_{L^2} $, on the availability of lower and upper estimates for $ \| D^{m-1} \uu(\cdot,t) \|_{L^2} $, and on monotonicity results for $ \| D^{m} \uu(\cdot,t) \| _{L^2}$, see Lemma \ref{monotonicity-lemma}. There are far fewer results for estimates on lower bounds for decay as in Theorems \ref{lower-bound-decay-first} and \ref{lower-bound-decay-second} (i.e. for derivatives of solutions $\uu$), see Deng and Shang \cite{MR4220767} and Oliver and Titi \cite{MR1749867} for the Navier-Stokes equations; Ahn \cite{MR4103923} and Ye, Jia and Dong \cite{MR4135099} for MHD equations and related models, and Bae and Biswas \cite{MR3457535} for rather general equations. 

\begin{Remark} \label{comparison} There are well-known methods, different from the one we use here, which have led to results as in Theorems \ref{upper-bound-decay}, \ref{lower-bound-decay-first} and \ref{lower-bound-decay-second} for specific cases. In order to illustrate the differences, we briefly describe here the Fourier Splitting Method and the Gevrey estimates and radius of analiticity method. In the Fourier Splitting Method, developed by M.E. Schonbek  \cite{MR571048},   \cite{MR775190},  \cite{MR837929}, a differential inequality is established in frequency space, bounding the $L^2$ or $\dot{H} ^m$ norm by the average of the solution, or its $(- \Delta)  ^{  \frac{m}{2} }$ derivative,   over small enough frequencies. These estimates are obtained for   solutions to a regularized equation and then proved to hold for weak solutions through a standard limiting process, see M.E. Schonbek \cite{MR1312701} and M.E. Schonbek and Wiegner \cite{MR1396285} for the case of the Navier-Stokes equations. A different method, involving Gevrey estimates, radius of analiticity of solutions and analytic function theory, was developed by Oliver and Titi \cite{MR1749867} to establish similar estimates also for the Navier-Stokes equations. These ideas were further developed and extended for other functions spaces by Bae and  Biswas \cite{MR3457535} and Biswas \cite{MR2964642}, to study equations with  fractional Laplacian and quadratic analytic nonlinearities, thus providing a more general context on which to apply  this method. Note that in both these methods the essential estimates are obtained in frequency space.
\end{Remark}

\begin{Remark}
\label{remark-decay-character}
We now briefly discuss hypotheses (\ref{eqn:h1}) and (\ref{eqn:h4}) concerning upper and lower bounds for the decay of $\Vert \uu (t)  \Vert _{L^2}$, which are  starting points for the bootstrap methods that are central to our proofs.

Many results concerning upper and lower bounds for $\Vert \uu (t)  \Vert _{L^2}$, where $\uu$  is a solution to an  equation as in (\ref{eqn:main-equation}), are known; for a rather complete description of such results for the Navier-Stokes equations, we refer to Brandolese and M.E. Schonbek \cite{MR3916783}. In order to have a finer understanding of the relation between decay rate and initial data, Bjorland and M.E. Schonbek \cite{MR2493562}  introduced the decay character $r^{\ast} = r^{\ast} (\uu _0)$, for $\uu_0 \in L^2 (\RR ^n)$, defined  as $r^{\ast} (\uu _0)$ being the unique  $r \in \left( -\frac{n}{2}, \infty \right)$ such that $0 < P_r (v_0) < \infty$, for 

\begin{displaymath}
P_r(\uu_0) = \lim _{\rho \to 0} \rho ^{-2r-n} \int _{B(\rho)} \bigl |\widehat{\uu_0} (\xi) \bigr|^2 \, d \xi,
\end{displaymath}
provided this limit exists. Here, $B(\rho)$ denotes the ball at the origin with radius $\rho$. Using the decay character, Bjorland and M.E. Schonbek were able to establish tight decay estimates for solutions to the heat equation and as an application they proved that for solutions to the Navier-Stokes equations, when $- \frac{n}{2} < r^{\ast} \leq 1$ we have that 
 
\begin{displaymath}
 C_1 (1 + t) ^{- \left( \frac{n}{2} + r^{\ast} \right)} \leq \Vert \uu (t) \Vert _{L^2} ^2 \leq C_2 (1 + t) ^{- \left( \frac{n}{2} + r^{\ast} \right)},
\end{displaymath}
and for $r^{\ast} \geq 1$, we have 

\begin{displaymath}
\Vert \uu (t) \Vert _{L^2} ^2 \leq C (1 + t) ^{- \left( \frac{n}{2} + 1 \right)}.
\end{displaymath}
Note that for $- \frac{n}{2} < r^{\ast} \leq 1$ these estimates imply that the decay rates in (\ref{eqn:h1}) and (\ref{eqn:h4}) are equal, so Theorem \ref{lower-bound-decay-first} applies. However, as there are no known results for lower bounds when $r ^{\ast} > 1$, we have no information on (\ref{eqn:h4}) and Theorem \ref{lower-bound-decay-second} may apply.

The decay character can be explicitly computed in many important cases, see C\'ardenas and Niche \cite{CARDENAS2022125548} and Ferreira, Niche and Planas \cite {MR3565380}. However, the limit defining $P_r (\uu _0)$ may not exist, as is required in the definition.  Brandolese \cite{MR3493117} constructed initial data in $L^2$, highly oscillating near the origin, for which the decay character cannot be defined. He also provided  generalizations of the decay character, which he called  upper and lower decay characters, which lead to upper and  lower estimates for decay rates of solutions to linear systems for diagonalizable operators, as defined in Niche and M.E. Schonbek \cite{MR3355116}. Moreover, he showed that the decay character exists \ if and only if $\uu_0$ belongs to a certain subset of a homogeneous Besov space, i.e. $\uu_0  \in  \dot{\mathcal{A}} ^{- \left(\frac{n}{2} + r^{\ast} \right)} _{2, \infty} \subset \dot{B} ^{- \left(\frac{n}{2} + r^{\ast} \right)} _{2, \infty}$. In that case solutions to the  linear system are such that

\begin{displaymath}
C_1 (1 + t)^{- \frac{1}{\a} \left( \frac{n}{2} + r^{\ast} \right)} \leq \Vert \uu (t) \Vert _{L^2} ^2 \leq C_2 (1 + t)^{- \frac{1}{\a} \left( \frac{n}{2} + r^{\ast} \right)},
\end{displaymath}
if and only if the decay character $r^{\ast} = r^{\ast} (\uu_0)$ exists. These estimates are needed for establishing decay rates for nonlinear equations as  in (\ref{eqn:main-equation}) and, as a consequence, exert influence on whether the rates on  (\ref{eqn:h1}) and (\ref{eqn:h4}) may be equal or not, thus leading the potential use of either    Theorem \ref{lower-bound-decay-first} or Theorem \ref{lower-bound-decay-second}. $\Box$
\end{Remark}

Our next results provide ``reverse''estimates: under rather mild hypotheses, decay of the norm of $D \uu$ leads to decay of the $L^2$ norm of the solution $\uu$.

\begin{Theorem}
\label{upper-to-lower-above}
Assume a solution $\uu$ is such that $\lim_{t \rightarrow \infty} \| \uu(t) \| _{L^2} =  0 $, that (\ref{eqn:h2}) holds with $m = 0$ and that

\begin{equation}
\label{eqn:upper-bound-derivative}
\| D \uu(t) \| _{L^2} \leq C_{1} \, t^{-\alpha_{1}}, \quad \forall t > T_{1}.
\end{equation}
Then:

\begin{enumerate}
\item if $\alpha_{1} >  \frac{1}{2}  $  and (\ref{eqn:h3}) holds  with $ m = 0 $ for some $ \beta > \alpha_{1} +  \frac{1}{2}  $,
then

\begin{displaymath}
\| \uu (t) \| _{L^2}  \leq C_{0} \,  t^{- \alpha_{1}  +  \frac{1}{2} }, \quad \forall t > T_{0}.
\end{displaymath}
The constant $C_0$ is such that $C_0 = C_0 \left( C_{1}, \alpha_{1}, \beta, \sigma_0, F_0, \| \uu(a_{\ast})\| _{L^2} \right)$, with $ a_{\ast} $  as in Remark \ref{rmk-a-star}; $T_0$ is such that $T_0 = T_0 \left( t_{\ast}, T_{1}, \alpha_{1}, \beta, \nu,\tau_{0}, \sigma_{0}, F_0, a_{\ast}, g_0 \right)$;

\item  if $\alpha_{1} >  \frac{1}{2}  $  and (\ref{eqn:h3}) holds  with $ m = 0 $ for  $ \beta = \alpha_{1} +  \frac{1}{2}  $,
then for every $\epsilon > 0$ we have

\begin{displaymath}
\| \uu(t) \| _{L^2} \leq     C_{0} t^{- \alpha_{1} +  \frac{1}{2}  + \epsilon}, \quad \forall t > T_{0},
\end{displaymath}
such that $C_0 = C_0 \left( \epsilon, C_{1}, \alpha_{1}, \beta, \sigma_0, F_0,  \| \uu(a_{\ast})\| _{L^2} \right)$, where $a_{\ast}$  is as in Remark \ref{rmk-a-star}; $T_0$ is such that $T_0 = T_0 \left( t_{\ast}, T_{1}, \alpha_{1}, \nu, \tau_{0}, \sigma_{0}, F_0, a_{\ast}, g_0 \right)$, but $T_0$ does not depend on $\epsilon$.

\end{enumerate}

\end{Theorem}

\begin{Theorem}
\label{upper-to-lower-below}
Assume a solution $\uu$ is such that $\lim_{t \rightarrow \infty} \| \uu(t) \| _{L^2} =  0 $, that (\ref{eqn:h2}) holds with $m = 0$, that (\ref{eqn:upper-bound-derivative}) holds for some $\alpha_1 >  \frac{1}{2} $ and that

\begin{equation}
\label{eqn:hypotheses-upper-to-lower-below}
\| D \uu(t) \| _{L^2} \geq D_{1} \,\nu^{-\frac{1}{2}}\, t^{-\alpha_{1}}, \quad \forall t > t_{1},
\end{equation}
for some $D_1, t_1 > 0$. If (\ref{eqn:h3}) holds  with $ m = 0 $ for some $ \beta > \alpha_{1} +  \frac{1}{2}  $, then

\begin{displaymath}
\| \uu(t) \| _{L^2}  \geq  D_0 \, t^{- \alpha_{1}  +  \frac{1}{2} }, \quad \forall  t > t_{0}.
\end{displaymath}
\mbox{ } 

We have that $D_0 = \frac{D_1}{2 \sqrt{\alpha_1 -  \frac{1}{2} }}$, while $t_0 = t_0 ( t_{\ast}, C_1, D_1, t_{1}, T_1, \alpha_{1}, \beta, \nu,\tau_{0}, \sigma_{0}, F_0$, $ a_{\ast}, ||\uu(a_{\ast})||_{L^2}, g_0) $,  where $a_{\ast}$  is as in Remark \ref{rmk-a-star}.
\end{Theorem}
Even though these results seem natural, we have not been able to find them in the literature. 

\begin{Remark}
\label{rmk-a-star}
It will become apparent in the proofs of Theorems \ref{upper-to-lower-above} and \ref{upper-to-lower-below} that we only need $\liminf_{t\,\rightarrow\,\infty} \| \uu(t) \|_{L^2} = 0$.  However, using a Gronwall-type argument it can be shown that, for some $a_{\ast} \gg 1 $,  (depending in general on $t_{\ast}, \tau_0, \nu$ and the function $g_0$), we have
\begin{displaymath}
\| \uu(t) \| _{L^2} \leq \| \uu(a) \|_{L^2} + \int_a ^t \| \ff (\tau) \| _{L^2} \,d \tau,  \quad
\forall  t > a \geq a_{\ast}.
\end{displaymath}
Since $ \beta > 1 $, it follows that  $ \Vert \ff (\cdot) \Vert  _{L^2} \in L^{1}(a, \infty) $;
this then leads to
\begin{displaymath}
\limsup_{t \rightarrow \infty} \| \uu(t)\ \| _{L^2}  \leq  \liminf_{t\,\rightarrow\,\infty} \| \uu(t) \|_{L^2}.
\end{displaymath}
Hence,  the limit  $ \lim_{t\,\rightarrow\,\infty} \| \uu(t) \| _{L^2} $  does exist.
\end{Remark}

\begin{Remark}
We note that when $\ff = 0$, the hypothesis on $\| \uu (t) \|_{L^2}$ in Theorem \ref{upper-to-lower-below}  can be dropped.
\end{Remark}

This article is organized as follows. In Section \ref{proofs} we prove Theorems  \ref{upper-bound-decay}, \ref{lower-bound-decay-first}, \ref{lower-bound-decay-second}, \ref{upper-to-lower-above}, and \ref{upper-to-lower-below}, as well as the useful Monotonicity Lemma \ref{monotonicity-lemma}. In Section \ref{examples} we apply our results  to obtain decay rates for advection-difussion equations and for MHD equations.

\section{Proofs}
\label{proofs}

{\bf Proof of Theorem \ref{upper-bound-decay}:} we follow the ideas in Hagstrom, Lorenz, J. Zingano and P. Zingano \cite{MR4021907}.  Let $ \gamma > 2 \alpha$ and $a \geq \max \{1, t_{\ast}, T_{0}, \tau_0,  \sigma_{0} \}$. Taking dot product  of (\ref{eqn:main-equation}) with $2 (t - a)^{\gamma}   \uu(x,t)$ and integrating on $\mathbb{R}^{n}  \times (a, t) $ we obtain, for all $t > a$

\begin{align*}
(t - a)^{\gamma} \, \| \uu(t) \| _{L^2} ^{2} & + 2 \nu \int_{a}^{t} (\tau - a)^{\gamma} \, \|D \uu(\tau) \| _{L^2} ^{2} d \tau = \gamma 
\int_{a}^{t} (\tau - a)^{\gamma - 1} \, \| \uu(\tau) \|^{2} _{L^2} d\tau \\ & + 2 \int_{a}^{t} (\tau - a)^{\gamma} \left( \int_{ \mathbb{R}^{n}} \langle \uu(x,\tau), \ff (x,\tau) - \mathbb{G}(\uu) \,\rangle dx \right)  d \tau \\ & \leq  \gamma \int_{a}^{t} (\tau - a)^{\gamma - 1} \, \| \uu(\tau) \|^{2} _{L^2}  d \tau \\ & +  2 \int_{a}^{t} (\tau - a)^{\gamma} \,  \left(  f (\tau)  \| \uu(\tau) \|_{L^2}  + g_0 (\tau) \| D \uu(\tau) \| _{L^2} ^{2} \right) d\tau \\ & \leq C \int_{a}^{t} (\tau - a)^{\gamma - 2 \alpha - 1} \, d\tau + 2 \int_{a}^{t} (\tau - a)^{\gamma} \, g_0 (\tau) \, \| D \uu(\tau) \|  _{L^2} ^{2} d \tau
\end{align*}
where we used (\ref{eqn:h1}), (\ref{eqn:h2}) and (\ref{eqn:h3}), where $f(\tau) = \Vert \ff (\tau) \Vert _{L^2}$ and $C = C_0 (\gamma C_0 + 2 F_0 )$. As $g_0 (\infty) = 0$ we then obtain, increasing $a$ if necessary, 

\begin{displaymath}
(t - a)^{\gamma} \, \| \uu(t) \| _{L^2} ^{2} + \nu \int_{a}^{t} (\tau - a)^{\gamma} \, \| D \uu(\tau) \|_{L^2} ^{2} d\tau \leq 
E_0 (t - a)^{\gamma - 2 \alpha}
\end{displaymath}
for all $t > a $,  where $ E_0 = \frac{C}{\gamma - 2 \alpha}$. This leads to 
\begin{equation}
\label{eqn:first-ineq-thm-upper-bound}
\int_{a}^{t} (\tau - a)^{\gamma} \, \| D \uu(\tau) \|_{L^2} ^{2} d \tau \leq  E_0 \,\nu^{- 1} \, (t - a)^{\gamma - 2 \alpha}
\end{equation}
for all $ t > a$, and any sufficiently large $ a \geq \max \{ 1, t_{\ast}, T_0, \tau_0, \sigma_0 \}$. Differentiating (\ref{eqn:main-equation}) with respect to $x_{\ell} $, multiplying by $2 (t - a)^{\gamma + 1}  D_\ell \uu(x,t)$, integrating on $\mathbb{R}^{n}  \times (a, t) $ and summing over $\ell$ we similarly obtain (increasing $a$ if necessary, so that  $ a \geq \max\,
\{1,  t_{\ast}, T_0, T_1, \tau_0, \tau_{1}, \sigma_0, \sigma_1, \} $), 

\begin{align}
\label{eqn:second-ineq-thm-upper-bound}
(t - a)^{\gamma + 1} \,
\| D \uu(t) \|_{L^2} ^{2} & +  \nu \int_{a}^{t} (\tau - a)^{\gamma + 1} \, \| D^{2} \uu(\tau) \|_{L^2} ^{2} d\tau \notag \\ & \leq (\gamma + 1)  \int_{a}^{t} (\tau - a)^{\gamma} \, \| D \uu(\tau) \|_{L^2} ^{2} \, d\tau  \notag \\ & +  2 \int_{a}^{t} (\tau - a)^{\gamma + 1} \, \|  D\uu(\tau) \|_{L^2}  \| D  \ff (\tau) \|_{L^2}  d\tau
\end{align}
for all $t > a$. If $ \beta = \alpha + 1 $,  we then have 
\begin{align}
\label{eqn:third-ineq-thm-upper-bound}
(t - a)^{\gamma + 1} \, \| D  \uu(t) \|_{L^2} ^{2} & + \nu \int_{a}^{t} (\tau - a)^{\gamma + 1} \, \| D^{2} \uu(\tau) \|_{L^2} ^{2} d\tau \notag \\ & \leq (\gamma  + 1 + \nu)  \int_{a}^{t} (\tau - a)^{\gamma} \, \| D \uu(\tau) \|_{L^2} ^{2} d \tau \notag \\ & +  \nu^{-  1}  \int_{ a}^{  t} (\tau - a)^{\gamma + 2} \, \|  D   \ff (\tau)    \|_{L^2}^{  2}
\, d\tau
\end{align}
while,  if $ \beta > \alpha + 1 $,  we obtain from (\ref{eqn:second-ineq-thm-upper-bound}) and (\ref{eqn:h3}) with $ m = 1 $ that
\begin{align}
(  t - a)^{\gamma + 1} \, \|  D   \uu(t)   \|_{L^2} ^{  2} & +  \nu \int_{   a}^{  t} (\tau - a)^{\gamma + 1} \, \|  D^{2}  \uu(\tau)   \|_{L^2} ^{  2} d\tau  \notag \\ &  \leq (\gamma + 2) \int_{   a}^{  t}   (\tau - a)^{\gamma} \, \|  D  \uu(\tau)   \|_{L^2} ^{  2} \, d\tau \notag \\ & +     \int_{   a}^{  t}   (\tau - a)^{\gamma + 2} \, \|  D     \ff (\tau)    \|_{L^2} ^{  2}
\, d\tau \notag \\ & \leq (\gamma + 2)    \int_{   a}^{  t}   (\tau - a)^{\gamma} \, \|  D  \uu(\tau)   \|_{L^2} ^{  2} \, d\tau  \notag \\ & + F_{1}^{  2} \, \frac{1}{a^{2 \delta}}    \int_{   a}^{  t}   (\tau - a)^{\gamma   -   2  \alpha   -   1} \, d\tau
\end{align}
where $ \delta = \beta - (\alpha + 1) $.  Therefore, in the case $ \beta > \alpha + 1 $ we obtain, using (\ref{eqn:first-ineq-thm-upper-bound}) and once again  increasing $a$ if necessary 
\begin{equation}
(  t - a)^{\gamma + 1} \, \|  D   \uu(t)   \|_{L^2} ^{  2} +  \nu \int_{   a}^{  t} (\tau - a)^{\gamma + 1} \, \|  D^{2}  \uu(\tau)   \|_{L^2} ^{  2} d\tau
 \leq  E_{1} \, \nu^{-  1} (  t - a)^{\gamma   -   2  \alpha}
\end{equation}
for all $       t    > a $, with $E_{1}       = (\gamma + 2)^{2}    C_{0}^{  2} $. We have thus obtained our estimate for $m  =1$ and also

\begin{equation}
\int_{   a}^{  t}   (\tau - a)^{\gamma + 1} \, \|  D^{2}  \uu(\tau)   \|_{L^2} ^{  2} d\tau  \leq 
E_{1} \, \nu^{-  2} (  t - a)^{\gamma   -   2  \alpha},
\end{equation}
for all $     t    > a $.  We can use this first estimate to bootstrap our computations for the case $m = 2$ and proceed this way up to $\widehat{m} = m$. If $ \beta = \alpha + 1 $, we proceed  similarly starting from (\ref{eqn:third-ineq-thm-upper-bound}) to obtain the estimate for $ m = 1 $, then bootstrapping. $\Box$

We prove now a Monotonicity Lemma which will be necessary in the proof of the next estimates. This Lemma extends Theorem B in Braz e Silva, J. Zingano e P. Zingano \cite{MR3907942} 

\begin{Lemma}[Monotonicty Lemma]
\label{monotonicity-lemma}
Assume (\ref{eqn:h1}), (\ref{eqn:h2}) and (\ref{eqn:h3}) hold and $ \beta \geq \alpha + 1$. Then,  for every $ 0 \leq m \leq \widehat{m}  $ we have that 

\begin{equation}
\frac{d}{  d t  } \left(  \|  D^{m}    \uu(t)   \| _{L^2} ^{  2} +\,  K_{m, \alpha,\beta}  \, t^{  -  \alpha   -  \beta   -  m   +  1} \right) 
    \leq  0, \quad   \,
\forall  \, t > a_{m}
\end{equation}
where $K_{   m, \alpha,\beta} =   \frac{2  C_{m}   F_{m}}{\alpha + \beta + m - 1} $, with $      a_{m}    > t_{\ast} $ such that $a_m = a_m \left(T_{m},    \tau_{m},    \sigma_{m},   \nu, g_m \right)$. The constants $C_m, F_{m}$ referred to here are those in (\ref{eqn:h1}), (\ref{eqn:h3}) and Theorem \ref{upper-bound-decay}.
\end{Lemma}

\begin{Remark}
If $ \ff = \mbox{\bf 0} $, it then follows that  $ \|  D^{m}      \uu(t)  \|_{L^2}  $ is monotonically decreasing in  the interval $ (a_{m}, \infty) $, since in this case  we have  $    K_{m}  = 0 $, as $F_{m}       = 0 $.  
\end{Remark}

{\bf Proof of Lemma \ref{monotonicity-lemma} :} Let $   a   =  \max  \{  t_{\ast},  T_0,    \tau_0,   \sigma_{0}      \} $, and let $  t > a $. 
From the energy identity 
\begin{displaymath}
\|\,\uu(t)   \| _{L^2} ^{  2} \,+     2    \nu   \int_{   a}^{  t}    \|  D  \uu(\tau)  \|_{L^2} ^{  2}
  d\tau  =  \|  \uu(a)  \|_{L^2} ^{  2} \,+     2   \int_{   a}^{  t}   \int_{\mathbb{R}^{n}}\!     \langle   \uu(x,\tau), \ff (x,\tau)    -\mathbb{G}(\uu) \,\rangle  dx\, d\tau
\end{displaymath}
and (\ref{eqn:h2}) we obtain 
\begin{align*}
\frac{d}{  d t  } \, \|\,\uu(t)   \|_{L^2} ^{  2} & \,+     2    \nu \, \|  D  \uu(t)  \|_{L^2} ^{  2}
 =   2  \int_{\mathbb{R}^{n}} \langle   \uu(x,t), \ff (x,t)    -\mathbb{G}(\uu) \,\rangle  dx \\ & \leq  
2\, \|\,\uu(t)  \|_{L^2}   \|   \ff (t)  \|_{L^2} \,+\, 2\, g_0   (t)   \|  D  \uu(t)  \|_{L^2} ^{  2}.
\end{align*}
Then, after increasing $  a  $ if necessary, and using (\ref{eqn:h1}) and (\ref{eqn:h3})  we have that 
\begin{displaymath}
\frac{d}{  d t  } \, \|\,\uu(t)   \|_{L^2} ^{  2} \,+     \nu  
\|  D  \uu(t)  \|_{L^2} ^{  2}  \leq   2  C_0   F_0   t^{  -  \alpha   -  \beta}
\end{displaymath}
for all $  t > a $. Setting $ a_0 = a $, this is the result for $  m = 0 $. For $ 1 \leq m \leq \widehat{m} $, we proceed in a similar way by  taking $   a =  \max  \{  t_{\ast},  T_m,   \tau_{m},   \sigma_{   m}   \}$. Then through (\ref{eqn:h2})  we get
\begin{align}
\frac{d}{  d t  } \, \|  D^{m}    \uu(t)   \|_{L^2} ^{  2} \,+     2    \nu \, \|  D^{m+1}  \uu(t)  \|_{L^2} ^{  2}  & \leq  2  \|  D^{m}  \uu(t)   \|_{L^2}      \|  D^{m}     \ff (t)   \|_{L^2}  \notag \\ & +\, 2\, g_m   (t)   \|  D^{m+1}  \uu(t)   \|_{L^2} ^{  2}
\end{align}
for all $ t    > a $. Increasing $     a     $ if necessary,  through (\ref{eqn:h3}) and Theorem \ref{upper-bound-decay} we obtain
\begin{displaymath}
\frac{d}{  d t  } \, \|  D^{m}   \uu(t)   \|_{L^2} ^{  2} \leq  2     
\|  D^{m}  \uu(t)   \|_{L^2}       \|  D^{m}     \ff (t) \|_{L^2}  \leq     2    C_{m}   F_{m}     
t^{  -  \alpha   -  \beta   -  m},  \quad   \, \forall  \, t    > a. 
\end{displaymath}
Setting $ a_m = a $, this gives the result for $ D^m \uu(t) $, completing the proof.  $\Box$

\begin{Remark} From now on we use the notation

\begin{equation}
\label{eqn:zm}
z_{m}(t)  =  \|  D^{m}    \uu(t)   \|_{L^2} ^{  2}
  +\, K_{m, \alpha,\beta}  t^{  -  \alpha   -  \beta   -  m   +  1}, \quad t > a_{m}. 
\end{equation}
If $      \beta \geq \alpha + 1 $, Lemma \ref{monotonicity-lemma} implies that $ z_{m}$ is smooth and monotonically decreasing in the interval $(a_{m}, \infty) $. 
\end{Remark}

{\bf Proof of Theorem \ref{lower-bound-decay-first}:}  Let $  m = 1 $ and $t  > t_{1}    =\,  \max\,\{   t_{\ast},   t_{0},   \tau_{0},   \sigma_0, a_{1},   T_0   \} $. Given $ T       =    M   t $,  where $ M   > 1 $ will be chosen later, we obtain  from (\ref{eqn:main-equation})
 \begin{displaymath}
\|\uu(T)   \|_{L^2} ^{  2} +     2   \nu \!  \int_{t}^{  T}\! \|  D  \uu(\tau)  \|_{L^2} ^{  2}  \, d\tau =  \|\uu(t)   \|_{L^2} ^{  2} +     2\!  \int_{t}^{  T}\!\!  \int_{\mathbb{R}^{n}} \langle \uu(x,\tau),  \ff (x,\tau) - \mathbb{G}(\uu) \rangle dx\,d\tau.
\end{displaymath}
Using (\ref{eqn:h2}) with $ m = 0 $, we obtain
\begin{align*}
\| \uu(T)   \|_{L^2} ^{  2} & +     2   \nu   \int_{t}^{  T}  \|  D \uu(\tau)  \|_{L^2} ^{  2}  \, d\tau \,\geq  \|\,\uu(t)   \|^{  2} \\ & - 2  \int_{t}^{  T} \left( \|  \uu(\tau)   \|_{L^2}  \|      \ff (\tau)  \|_{L^2}  +     g_{0}   (\tau) \|  D \uu(\tau)   \|_{L^2} ^{  2} \right)  \,d\tau.
\end{align*}
Increasing $t_{1}$ if necessary, so that $ g_{0}      (\tau) < \nu      $ for all $ \tau > t_{1} $, through (\ref{eqn:h1}), (\ref{eqn:h3})  and (\ref{eqn:h4})  we obtain  
\begin{align*}
4  \nu   \int_{t}^{  T} \|  D  \uu(\tau)  \|_{L^2} ^{  2}  \, d\tau & \geq  \|\,\uu(t)   \|_{L^2} ^{  2} \,- \|\,\uu(T)   \|_{L^2} ^{  2} \,-\,2   \int_{t}^{  T}  \|\,\uu(\tau)   \|_{L^2}  \|    \ff (\tau)  \|_{L^2}  \, d\tau \\ & \geq  \left(     D_0 ^{2}  -   C_{0}^{  2}   M^{-  2  \alpha} -   \kappa_0 \, t^{  -  \delta}  \right) \, t^{  -  2  \alpha}
\end{align*}
where  $     \kappa_0    =   \frac{2\,C_0  F_0}{\alpha + \beta - 1}$ and $      \delta = \beta - \alpha - 1 > 0$. Choosing $ M       =  \bigl(  \frac{2 C_0}{D_0}  \bigr)^{ 1/\alpha} $ and again increasing $ t_{1}$ if necessary, so that  $ \kappa_0 \, t_{1}^{  -  \delta}  \leq  \frac{1}{4} D_0 ^{2}$, by Lemma \ref{monotonicity-lemma} we obtain 
\begin{displaymath}
4\,\nu\,T      z_{1}(t)  \geq  4   \nu   \int_{t}^{  T}  z_{1}(\tau) \, d\tau  \geq  4   \nu   \int_{t}^{  T} \|  D  \uu(\tau)  \|_{L^2} ^{  2} \, d\tau  \geq    \frac{1}{2} D_0 ^2 \, t^{  -  2  \alpha}.
\end{displaymath}
Therefore if $     t       > t_{1}$ 
\begin{displaymath}
z_{1}(t) \geq \frac{D_0 ^2}{8 \nu}  \, T^{  -  1}   t^{  -  2  \alpha}  =  \frac{D_0  ^2}{8     M     \nu } \,  t^{  -  2  \alpha   -  1},
\end{displaymath}
or equivalently
\begin{displaymath}
\|  D  \uu(t)   \| _{L^2} ^{  2} \,\geq \left(  \frac{D_0  ^2}{8     M     \nu }   -   K_{1, \alpha, \beta}      t^{  -  \delta} \right) \, t^{  -  2  \alpha   -  1}
\end{displaymath}
for all $   t    > t_{1} $. Increasing $t_{1}$ if necessary, we obtain
\begin{displaymath}
\|  D  \uu(t)  \| _{L^2}  \geq  D_1  \nu ^{- \frac{1}{2}} t^{  -  \alpha   -   \frac{1}{2} },
\end{displaymath}
for $  t > t_{1} $ and $D_1 = \frac{D_0}{3 \sqrt{M}}$.  We have thus proved our result for $ m = 1 $.   If $     \widehat{m} = 1 $  we are done; otherwise, we proceed with $ m = 2 $ in a similar way. Let $t >       t_{2}    =  \max\,\{   t_{1},       \tau_{1},       \sigma_{1},      a_{2}, T_1 \} $  and $      T       =       M     t $, where $ M       = \left( \frac{2  C_1}{D_1} \right)^{1/\alpha}     $. Increasing  $     t_{2}    $ if necessary, we then get from (\ref{eqn:h2}), (\ref{eqn:h3}), Theorem \ref{upper-bound-decay} and the estimate for $m  =1$
\begin{align*}
4 \nu\!\!   \int_{t}^{  T}        \|  D^{2}  \uu(\tau)  \|_{L^2} ^{  2}  \, d\tau & \geq  \|  D  \uu(t)   \|_{L^2} ^{  2} \,-     \|  D  \uu(T)   \|_{L^2} ^{  2} \,-\,2\!   \int_{t}^{  T}\!  \|  D  \uu(\tau)   \|_{L^2}   \|  D            \, \ff (\tau)  \|_{L^2}  \, d\tau \\ & \geq   \left(\frac{D_1 ^2}{\nu}  -   \frac{C_{1}^{  2}   M^{-  2  \alpha  -  1}}{\nu}  -   \frac{\kappa_1}{\nu^{ \frac{1}{2}}  }      t^{  -  \delta} \right)      t^{  -  2  \alpha   -  1} \geq \frac{1}{2} D_1 ^{2} \nu ^{-1} \, t^{  -  2  \alpha   -  1}.
\end{align*}
Now, for $z_2$ as in Lemma \ref{monotonicity-lemma}, we have that $     z_{2}(t)       \geq  \|  D^{2}     \uu(t)  \|_{L^2} ^{  2}    $. We start with this bound and repeat bootstraping until $m = \widehat{m}$. $\Box$

{\bf Proof of Theorem \ref{lower-bound-decay-second} :} Let $      q = \frac{\eta}{\alpha} $ and   $   t_{1}    =\,  \max\,\{   t_{\ast},   t_{0}, 
  \tau_{0},   \sigma_0,    a_{1},  T_0   \} $, where these constants are as in (\ref{eqn:h1}), (\ref{eqn:h2}), (\ref{eqn:h3}), (\ref{eqn:h4}) and (\ref{eqn:zm}). For $ t > t_{1} $ we take $ T       =    M   t^{  q}    $,  where $ M   > 1 $ will be chosen later, so that we have 
\begin{displaymath}
\|\,\uu(\mbox{\small $T$})   \|_{L^2} ^{  2} \,+     2 \,\nu\!\!\!   \int_{t}^{  T} \! \|  D  \uu(\tau)  \|_{L^2} ^{  2}  \, d\tau    = 
\|\,\uu(t)   \|_{L^2} ^{  2} \,+     2\!\!   \int_{t}^{  T}\!\!\!\int_{\mathbb{R}^{n}}\! \langle \uu(x,\tau),     \ff (x,\tau) -    \mathbb{G}(\uu) \,\rangle dx\,d\tau.
\end{displaymath}
By (\ref{eqn:h2}) with $ m = 0 $ we obtain 
\begin{align}
\|\,\uu(T))   \|_{L^2} ^{  2} & +     2   \nu   \int_{t}^{  T}  \|  D \uu(\tau)  \|_{L^2} ^{  2}  \, d\tau \,\geq \|\,\uu(t)   \|_{L^2} ^{  2} \notag \\ & -     2   \int_{t}^{  T} \left(  \|  \uu(\tau)   \|_{L^2}  \| \ff (\tau)  \|_{L^2}  +    g_{0}   (\tau) \, \|  D \uu(\tau)   \|_{L^2} ^{  2} \right) \,d\tau.
\end{align}
Increasing $  t_{1}$, if needed,  so that $ g_{0}   (\tau) < \nu      $ for all $ \tau       > t_{1}   $, we have from  (\ref{eqn:h1}), (\ref{eqn:h3}) and  (\ref{eqn:h4}) that
\begin{align*}
4   \nu   \int_{t}^{  T}  \|  D  \uu(\tau)  \|^{  2}  \, d\tau & \geq  \|\,\uu(t)   \|^{  2}  \,-      \|\,\uu(T)   \|^{  2} \,-\,2  \int_{t}^{  T}  \|\,\uu(\tau)   \| \|      \ff (\tau)  \| \, d\tau \\ & \geq  \left(      D_0 ^{2}  -   C_{0}^{  2}   M^{-  2  \alpha}  -   \kappa_0 \, t^{  -  \delta_{1}} \right) \, t^{  -  2  \eta}
\end{align*}
where  $     \kappa_0    =   \frac{2\,C_0  F_0}{\alpha + \beta - 1}$ and $      \delta_{1}     = \beta -  (     2  \eta - \alpha + 1)    > 0$.  Choosing $ M       =  \left( \frac{2  C_0}{D_0} \right) ^{ \frac{1}{\alpha}} $ and increasing $ t_{1}  $ (if necessary) so that  $ \kappa_0 \, t_{1}^{  -  \delta_{1}}  \leq  \frac{D_0^{2}}{4}$, we obtain through Lemma \ref{monotonicity-lemma} 
\begin{displaymath} 
4\,\nu\,T      z_{1}(t)   \geq     4  \nu   \int_{t}^{  T}  z_{1}(\tau)  \, d\tau   \geq     4   \nu   \int_{t}^{  T} \|  D  \uu(\tau)  \|_{L^2} ^{  2} \, d\tau \geq    c^{  2}\, t^{  -  2  \eta}
\end{displaymath}
where  $     c^{  2}    =  \frac{D_0^{2}}{2}$. Therefore, $ z_{1}(t) \geq    \frac{ c^{  2}}{ 4 \nu T }   t^{  -  2  \eta}       = \frac{c^{  2}}{4  M \nu}     t^{  -  2  \eta   -  q} $ if $     t       > t_{1}$, or 
\begin{displaymath}
\|  D  \uu(t)   \|_{L^2} ^{  2} \,\geq  \left(    \frac{c^{  2}}{  4  M  }    -   K_{1, \alpha, \beta} \, \nu t^{  -  \varepsilon_{1}}  \right) \, \nu^{  -  1}      t^{  -  2  \eta   -  q}
\end{displaymath}
for all $   t    > t_{1} $, where $      \varepsilon_{1}    =      \beta - (     2  \eta - \alpha + q     )    > 0 $. Hence, increasing $     t_{1}   $ if necessary, we have  $ \|  D  \uu(t)  \|_{L^2}  \geq   D_1 \, \nu^{-   \frac{1}{2} } \, t^{  -  \eta   -  \frac{q}{2}} $ for $  t > t_{1} $, with  $      D_1 = \frac{D_0}{3 \sqrt{M}}$.  This shows the result for $     m = 1 $.  If $     \widehat{m} = 1 $,  the proof of this result  is complete; otherwise, we proceed with $ m = 2 $ in a similar way. Setting $      t_{2}    =  \max\,\{   t_{1},       \tau_{1},       \sigma_{1},      a_{2}, T_{ 1} \} $, let $      t       >    t_{2} $ and $    T          =       M t^{  q}     $, where $ M       = \left(\frac{2 C_1}{D_1} \right)^{1/\alpha}    $. As before, increasing $     t_{2}    $ if necessary, we then obtain from (\ref{eqn:h2}) with $m = 1 $, 
\begin{align*}
4\nu\!\!\int_{t}^{  T}       \|  D^{2}  \uu(\tau)  \|_{L^2} ^{  2}  \, d\tau    & \geq  \|  D  \uu(t)   \|_{L^2} ^{  2} \,-      \|  D  \uu(T)   \|_{L^2} ^{  2} \,-\,2\!\!   \int_{t}^{  T}       \|  D  \uu(\tau)   \|_{L^2}  \|  D  \, \ff (\tau)  \|_{L^2}  \, d\tau \\ &  \geq  \left(     D_1^{2}  -   C_{1}^{  2}   M^{-  2  \alpha  -  1} -   \kappa_1 \, \nu^{ \frac{1}{2} }     t^{  -  \delta_{  2}} \right)  \nu^{-  1}  t^{  -  2  \eta   -  q}
\end{align*}
by (\ref{eqn:h3}) and Theorem \ref{upper-bound-decay}  (both applied with $m = 1$), where $\kappa_{1}  = \frac{     2  C_{1}     F_{1}}{\alpha + \beta} $ and $      \delta_{  2}        = \beta    -    (     2  \eta - \alpha + q     ) > 0 $.  Increasing $     t_2 $ if necessary so that we have  $       \kappa_{1} \, \nu^{ \frac{1}{2} } \, t_{2}^{  -  \delta_{2}} <   \frac{ D_1 ^{2}}{4}$, using Lemma \ref{monotonicity-lemma}  we obtain 
\begin{displaymath}
4\,\nu\,T      z_{2}(t)  \geq     4  \nu   \int_{t}^{  T}  z_{2}(\tau)  \, d\tau   \geq    4   \nu   \int_{t}^{  T}  \|   D^{2}  \uu(\tau)  \|_{L^2} ^{  2}  \, d\tau  \geq     c^{  2}\,\nu^{-  1} \, t^{  -  2  \eta   -   q}
\end{displaymath}
for all $     t    > t_{2}    $, where  $     c^{  2}    =  \frac{D_1 ^2}{2}$. Therefore, $ z_{2}(t) \geq     \frac{c^{  2}}{4  \nu^{2} T} t^{  -  2  \eta   -  q} $,  or 
\begin{displaymath}
\|  D^{2}  \uu(t)   \| _{L^2} ^{  2} \,\geq \left(  \frac{c^{  2}}{  4  M  }   -   K_{2, \alpha, \beta} \, \nu^{  2}  \,t^{  -  \varepsilon_{2}}  \right) \, \nu^{  -  2}      t^{  -  2  \eta \,-\, 2  q}
\end{displaymath}
for all $   t    > t_{2} $, where $      \varepsilon_{2} =     \beta - (     2  \eta - \alpha + 2  (q - 1) + 1     )    > 0 $. Thus, increasing  $     t_{2} $ if needed,  we have $ \|  D^{2}  \uu(t)  \|_{L^2}  \geq    \frac{D_2}{\nu}  t^{  -  \eta \,-\, q} $ for all $      t    > t_{2} $, with  $      D_2 = \frac{D_1}{3 \sqrt{M}}$.  This completes  the proof  if $     \widehat{m} = 2 $. Otherwise we continue to the next level by considering the  energy estimate
\begin{align*}
\|  D^{2}  \uu(T)   \|_{L^2} ^{  2} & \,+     2   \nu  \int_{t}^{  T} \|  D^{3}  \uu(\tau)  \|_{L^2} ^{  2} 
\, d\tau    =  \|  D^{2}  \uu(t)   \|_{L^2} ^{  2}  \\ & + 2   \sum_{j\,=\,1}^{n} \sum_{\ell\,=\,1}^{n} \int_{t}^{  T}  \int_{\mathbb{R}^{n}} \langle  D_{   j}      D_{\ell} \uu(x,\tau),    D_{   j}      D_{\ell} \ff (x,\tau) -   D_{   j}      D_{\ell} \mathbb{G}(\uu) \,\rangle dx\,d\tau
\end{align*}
for $    t > t_{3}    =  \max\,\{  t_{2},       \tau_{2},       \sigma_{2}, a_{3}, T_2 \} $, where $      T       =       M     t^{  q}     $, $ M       = \left(     \frac{2  C_{2}}{D_2} \right)^{   1/\alpha}    $, repeating the steps above. Because of the condition on $ \beta $ in the hypotheses,  we can proceed in this way up to a last level, given by the energy estimate  for $       m = \widehat{m}    -    1 $
\begin{align*}
\| &  D^{m}  \uu(T)   \|_{L^2} ^{  2}  \,+     2   \nu   \int_{t}^{  T} \|  D^{m + 1}  \uu(\tau)  \|_{L^2} ^{  2}  \, d\tau =  \|  D^{m}  \uu(t)   \|_{L^2} ^{  2}  \\ & + 2    \sum_{\ell_{1}, \cdots, \ell_{m} =\,1}^{n}  \int_{t}^{  T}  \int_{\mathbb{R}^{n}} \langle D_{\ell_{1}} \cdots D_{\ell_{m}} \uu(x,\tau),    D_{\ell_{1}} \, \cdots D_{\ell_{m}} \left( \ff (x,\tau)    - \mathbb{G}(\uu) \right)  \,\rangle dx\,d\tau
\end{align*}
for $      t       >    t_{\widehat{m}}       = t_{m+1} =  \max\,\{   t_{m},       \tau_{m},       \sigma_{m},  a_{m+1}, T_m \} $, where $      T       =       M     t^{  q}     $, $ M = \left( \frac{2  C_{m}}{D_m } \right) ^{   1/\alpha}$.  Using (\ref{eqn:h2}), (\ref{eqn:h3}), Theorem \ref{upper-bound-decay} and Lemma \ref{monotonicity-lemma} we then obtain 
\begin{displaymath}
4 \nu\,T      z_{\widehat{m}}(t)  \geq     4   \nu   \int_{t}^{  T}  \|  D^{\widehat{m}}  \uu(\tau)  \|_{L^2} ^{  2}  \, d\tau \geq   \left(  D_m ^{2} -  C_{m}^{  2}       M^{  -  2  \alpha}    - \kappa_{m}    \nu^{   \frac{m}{2} }  \,t^{  -  \delta_{\widehat{m}}}  \right)  \nu^{-  m} \, t^{  -  2  \eta \,-\, m  q}
\end{displaymath}
where  $      \delta_{\widehat{m}}       = \beta - \left(  2  \eta - \alpha +  (     q-1)      m + 1  \right) > 0 $. Increasing $     t_{\widehat{m}}    $ if necessary, we obtain $ z_{\widehat{m}}(t) \geq    \frac{D_m ^2 }{8 M} \nu^{-\widehat{m}} \,  t^{  -  2  \eta \,-\, \widehat{m}  q} $ for all $       t    > t_{\widehat{m}} $,  or, in terms of $ \|  D^{  \widehat{m}}  \uu(t)  \|_{L^2}  $, by (\ref{eqn:zm})
\begin{displaymath}
\|  D^{  \widehat{m}}    \uu(t)  \|_{L^2} ^{  2}   \geq\, \left(  \frac{D_m ^2}{8M }   -  K_{\widehat{m},\alpha,\beta} t^{  -  \varepsilon_{\widehat{m}}}  \right)  \, \nu^{-  \widehat{m}} t^{  -  2  \eta \,-\, \widehat{m}  q}
\end{displaymath}
for all $      t      > t_{\widehat{m}} $, where $   \varepsilon_{\widehat{m}}    =   \beta - (  2  \eta - \alpha +  (     q-1)  \widehat{m}      + 1) $. As $       \varepsilon_{\widehat{m}}       > 0 $, the result follows.  $\Box$

{\bf Proof of Theorem \ref{upper-to-lower-above}:}  Let $      T_0 =\, \max\,\{  1,   t_{\ast}, \tau_0,   \sigma_0, T_{1}   \} $, and let  $ T       > t > T_{0} $. From (\ref{eqn:main-equation}), we obtain 
\begin{displaymath}
\|\,\uu(T)  \|_{L^2} ^{  2} \,+ 2 \nu \int_{t}^{  T} \|  D  \uu(\tau)  \|_{L^2}^{  2} d\tau \,=  \|\,\uu(t)  \|_{L^2} ^{  2} +     2   \int_{t}^{  T} \langle     \uu(x,\tau),  \ff (x,\tau)    -     \mathbb{G}(\uu)  \,\rangle dx\, d\tau,
\end{displaymath}
hence by (\ref{eqn:h2}) we have (increasing $ T_0    $  if necessary)
\begin{displaymath}
\|\,\uu(T)  \|_{L^2} ^{  2} \,+     4  \nu         \int_{t}^{  T} \|  D  \uu(\tau)  \|_{L^2} ^{  2} d\tau \,\geq  \|\,\uu(t)  \|_{L^2} ^{  2}  -    2  \int_{t}^{  T} \|\,\uu(\tau)   \|_{L^2} \|     \ff (\tau)   \|_{L^2}  d\tau.
\end{displaymath}
Using (\ref{eqn:upper-bound-derivative}) and letting $     T \rightarrow \infty $, as $ { \|\,\uu(T)  \|_{L^2}     \rightarrow 0} $ we obtain
\begin{displaymath}
\|\,\uu(t)  \|_{L^2} ^{  2} \leq  \frac{4  C_{1}^{  2}}{\,2  \alpha_{1}    - 1  }      t^{  1   -  2  \alpha_{1}}   +     2   \int_{t}^{  \infty} \|\,\uu(\tau)   \|_{L^2} \|      \ff (\tau)   \|_{L^2}  d\tau
\end{displaymath}
for all $t    > T_{0} $. If $\ff = 0 $, we are done. Otherwise, from the assumption on $\ff$, we have 
\begin{equation}
\label{eqn:first-estimate-theorem-upper-lower-above}
\|\,\uu(t)  \|_{L^2} ^{  2} \leq  \frac{4  C_{1}^{  2}}{\,2  \alpha_{1}    - 1  }      t^{  1   -\, 2  \alpha_{1}}  +    2      F_{0} \int_{t} ^{  \infty} \tau^{-  \beta} \|\,\uu(\tau)   \|_{L^2}  d\tau
\end{equation}
if $  t    > T_0 $.  By Remark \ref{rmk-a-star},  we have, for some $    a_{\ast}     \gg 1 $, depending on $    t_{\ast}, \tau_0,   \nu, g_0  $, that
\begin{equation}
\label{eqn:second-estimate-theorem-upper-lower-above}
\|\,\uu(t)   \|_{L^2} \,\leq\, M_0   = \|\,\uu(a_{\ast})  \|_{L^2} \,+     \int_{a_{\ast} }^{  \infty} \|   \ff (\tau)  \|_{L^2} \,d\tau \quad   \, \forall  \, t > a_{\ast}.
\end{equation}
Redefining $    T_0$ to be $    T_0 =\,\max\,\{  1,   t_{\ast}, \tau_0,   \sigma_0, T_{1},   a_{\ast}   \} $, from (\ref{eqn:first-estimate-theorem-upper-lower-above}), (\ref{eqn:second-estimate-theorem-upper-lower-above})  and (\ref{eqn:h3}) we obtain
\begin{displaymath}
\|\,\uu(t)  \|_{L^2} ^{  2}    \leq  \frac{4  C_{1}^{  2}}{\,2  \alpha_{1} - 1  }      t^{  1   -\, 2  \alpha_{1}}   + \frac{  2  F_0  M_0  }{\beta - 1}   t^{  1   -   \beta}
\end{displaymath}
for all $      t       > T_0 $.  If $\frac{ \beta    - 1}{2}   \geq   \alpha_{1}    -  \frac{1}{2} $, we are done;  otherwise, we obtain 
\begin{equation}
\label{eqn:third-estimate-theorem-upper-lower-above}
\|\,\uu(t)  \|_{L^2}    \leq\, M_{1}   t^{\,q_{\mbox{}_{1}}    (1 -  \beta  )},  \quad   \, \forall  \, t > T_{0}
\end{equation}
with $     q_{1}       =  \frac{1}{2}    $ and  $ M_{1}    = \sqrt{\frac{4  C_{1}^{  2}}{2  \alpha_1    - 1} + \frac{2  F_0  M_0}{\beta - 1}}$. Let $      \theta = \frac{\alpha_{1}   -  \frac{1}{2}}{\beta - 1}$. If 
\begin{displaymath}
\|\,\uu(t)  \| _{L^2}    \leq\, M_{k} \, t^{  q_{\mbox{}_{k}}    (1 -  \beta  )}, \quad   \, \forall  \, t > T_{0}
\end{displaymath}
for some  $    q_{k} \in [  0,   \theta) $, we obtain from  (\ref{eqn:first-estimate-theorem-upper-lower-above}) 
\begin{equation}
\label{eqn:fourth-estimate-theorem-upper-lower-above}
\|\,\uu(t)  \|^{  2} \leq \frac{4  C_{1}^{  2}}{\,2  \alpha_{1}    - 1  }      t^{  1   -\, 2  \alpha_{1}} +     \frac{  2  F_0  M_k  } {  (1    + q_{  k})  (\beta    - 1)  } \,   t^{  (1   +   q_{  k})  (1   -   \beta)},
\end{equation}
for all $      t       >    T_0 $. Let $q_{k+1}    =  \frac{1 + q_{k}}{2} $. If $ q_{k+1}  (\beta    -    1)  \geq      \alpha_{1}    -     \frac{1}{2}  $, we are done; if not,  we have 
\begin{displaymath}
\|\,\uu(t)  \| _{L^2}  \leq\, M_{k  +1} \, t^{  q_{k+1}}    (1 -  \beta  ), \quad   \, \forall  \, t > T_{0}
\end{displaymath}
where  $ M_{k  +1}    = \sqrt{\frac{4  C_{1}^{  2}}{2  \alpha_1    - 1} + \frac{F_0  M_k}{q_{k+1} (\beta - 1)}}$. The constants $q_{k}$  are given recursively by 
\begin{displaymath}
q_{k+1} \,=  \frac{\,1\,}{2} \, (  1 +   q_{k}), \quad    \, k   \geq   0,
\end{displaymath}
with $q_{0}     = 0 $, so that $ q_{k}      = 1    - 2^{  -  k}    $. Now,  in the case ({\em i\/}), we have $ \theta < 1 $ (since $ \beta    > \alpha_{1}    +  \frac{1}{2}  $), and there will  exist $ k_{\ast}    \geq 0   $ such that $ q_{k_{\ast}}  (\beta    -    1)  <      \alpha_{1}    -     \frac{1}{2}    $ and  $      q_{k_{\ast} +  1}  (\beta    -    1)  \geq      \alpha_{1}    -     \frac{1}{2}  $. For such $k$, (\ref{eqn:fourth-estimate-theorem-upper-lower-above}) leads to  
\begin{displaymath}
\|\,\uu(t)   \|_{L^2}  \leq\, M_{k_{\ast} +  1}     t^{  -  \alpha_{1}   +    \frac{1}{2} }
\end{displaymath}
for all $  t    > T_{0} $,. Finally,  in the case ({\em ii\/}), where $ \theta = 1 $, as  $ 2   \alpha_1    -    1 >  (1 + q_{k}   )  (\beta    -    1) $ for all $k$, the second term on the right hand side of (\ref{eqn:fourth-estimate-theorem-upper-lower-above}) will always decay slower than the first term. However, given $   \epsilon > 0 $, we have $ (1 + q_{k}   )  (\beta    -    1) > 2  \alpha_{1}    -   1 -  2  \epsilon $ for large $k$, and so  the bootstrap iteration  can stop there. $\Box$

\begin{Remark} For the sake of completeness, we prove here the estimate in Remark \ref{rmk-a-star}, which was used in the proof above.  From (\ref{eqn:h2}) and (\ref{eqn:h3}) with $\widehat{m} = 0$, let $     a_{\ast}  =\,
\max\,\{   t_{\ast}, \tau_0      \} $. Given $     t    > a > a_{\ast} $ we have 

\begin{displaymath}
\|\,\uu(t)  \|_{L^2} ^{  2}  +     2   \nu      \int_{a}^{  t} \|  D  \uu(\tau)  \|_{L^2} ^{  2} d\tau \,=     \|\,\uu(a)  \|_{L^2} ^{  2} + 2   \int_{a}^{ t} \int_{\mathbb{R}^{n}} \langle \uu(x,\tau),    \ff (x,\tau)    - \mathbb{G}(u) \,\rangle    dx\,d\tau
\end{displaymath}
so that, increasing $  a_{\ast}  $ (if necessary) so as to have $ g_{0}   (\tau) \leq \nu$ for all $ \tau    >    a_{\ast} $, we obtain 
\begin{displaymath}
\|\,\uu(t)  \|_{L^2} ^{  2} \,\leq     \|\,\uu(a)  \|_{L^2} ^{  2} + 2   \int_{a}^{  t} \|\,\uu(\tau)  \|_{L^2}  \|   \ff (\tau)  \|_{L^2} d\tau
\end{displaymath}
or,  in terms of $      \mbox{v}(t) = \|\,\uu(t)  \|_{L^2}^{  2} $      
\begin{displaymath}
\mbox{v}(t)     =     \mbox{v}(a)  \,+    2   \int_{a}^{  t} \mbox{v}(\tau)^{ \frac{1}{2} } \|   \ff (\tau)  \| d\tau.
\end{displaymath}
Now,  let $ \mbox{w}       \in    C^{1} ([  a, \infty)) $ be given by $ \mbox{w}^{  \prime}(t) =  2    \mbox{w}(t)^{ \frac{1}{2} }   
\|   \ff (t)  \|_{L^2} $ for $       t    > a $, and $ \mbox{w}(a)    =      \mbox{v}(a) $, that is
\begin{displaymath}
\mbox{w}(t)^{ \frac{1}{2} }  \,=     \|\,\uu(a)  \|_{L^2} \,+ \int_{a}^{  t} \|   \ff (\tau)  \|_{L^2} d\tau.
\end{displaymath}
Since $   \|\,\uu(t)  \|_{L^2} =   \mbox{v}(t)^{ \frac{1}{2} } \leq   \mbox{w}(t)^{ \frac{1}{2} }$,  we have proved the result. $\Box$\end{Remark}

{\bf Proof of Theorem \ref{upper-to-lower-below}:}  Let $      t_0 =\, \max\,\{  1,   t_{\ast}, \tau_0,   \sigma_0,
     t_{1}, T_{0}   \} $, where $ T_0 $ is as in Theorem \ref{upper-to-lower-above} (1), and let $ T       > t > t_{0} $.
From (\ref{eqn:main-equation}) we have 
\begin{displaymath}
\|\,\uu(T)  \|_{L^2} ^{  2} \,+     2  \nu         \int_{t}^{  T} \|  D  \uu(\tau)  \|_{L^2} ^{  2} d\tau \,=  \|\,\uu(t)  \|_{L^2} ^{  2} + 2   \int_{t}^{  T} \langle     \uu(x,\tau),       \ff (x,\tau)    -    \mathbb{G}(\uu) \,\rangle dx\, d\tau,
\end{displaymath}
so from \ref{eqn:h2}), we obtain, increasing $t_0$ if necessary 
\begin{displaymath}
\|\,\uu(T)  \|_{L^2} ^{  2} \,+     \nu \int_{t}^{  T} \|  D  \uu(\tau)  \|_{L^2} ^{  2} d\tau \,\leq  \|\,\uu(t)  \|_{L^2} ^{  2}  +    2  
\int_{t}^{  T} \|\,\uu(\tau)   \|_{L^2} \|      \ff (\tau)   \|_{L^2}  d\tau.
\end{displaymath}
Hence,  we have 
\begin{displaymath}
\|\,\uu(t)  \|_{L^2} ^{  2} \geq  \nu  \int_{t}^{  T} \|  D  \uu(\tau)  \|_{L^2} ^{  2} d\tau - 2   \int_{t}^{  T} \|\,\uu(\tau)   \|_{L^2}  \|      \ff (\tau)   \|_{L^2} d\tau
\end{displaymath}
so that, letting $   T    \rightarrow    \infty $
\begin{displaymath}
\|\,\uu(t)  \|_{L^2} ^{  2} \geq  \nu \int_{t}^{  \infty} \|  D  \uu(\tau)  \|_{L^2} ^{  2} d\tau -    2   \int_{t}^{  \infty} \|\,\uu(\tau)   \|_{L^2} \|      \ff (\tau)   \|_{L^2}  d\tau
\end{displaymath}
for all $   t    > t_0 $. From Theorem \ref{upper-to-lower-above} (item (1)), (\ref{eqn:hypotheses-upper-to-lower-below}) and (\ref{eqn:h3})  with $ m = 0 $, we obtain 
\begin{displaymath}
\|\,\uu(t)  \|_{L^2} ^{  2} \geq  \frac{ D_1 ^{  2}} {\,2  \alpha_{1}    - 1  }   t^{  1   -  2  \alpha_{1}} \,-  \frac{2  C_0      F_0}
{   \beta + \alpha_{1}  - 3/2  }   t^{  \frac{3}{2} \,-   \beta   -  \alpha_{1}} \geq  \frac{1}{\,2  \alpha_{1}  - 1  } \left(  D_1 ^{2} -  
2  C_{0}       F_0 \, t^{  -  \delta} \right)  t^{  1   -  2  \alpha_{1}}
\end{displaymath}
for all $       t    > t_{0} $, where  $      \delta =       \beta - (\alpha_{1}    +     \frac{1}{2} ) $.  Since $ \delta > 0 $, increasing 
$      t_0        $ (if necessary) we obtain 
\begin{displaymath}
\|\,\uu(t)  \|_{L^2} ^{  2} \geq  \frac{D_1 ^{  2}} {  4 \left( \alpha_{1}    -  \frac{1}{2} \right)   }   t^{  1   -  2  \alpha_{1}}
\end{displaymath}
for all $   t    > t_0 $. Thus, we have proved our result. $\Box$

\section{Examples}
\label{examples}

\subsection{Advection-diffusion equations} Given a smooth function $\bb = (b_1, b_2, \cdots, b_n)$ in $C^{\infty}(\mathbb{R}) $ and an arbitrary initial datum  $u_0  \in L^{2}  $,  let  
\begin{displaymath}
u \in C([0, \infty), L^{2} )  \cap  L^{\infty}_{loc}( (0, \infty), L^{\infty}(\mathbb{R}^{n}))
\end{displaymath}
be a solution to the problem
\begin{equation}
\label{eqn:adv-diff} 
u_t + \bb (u) \cdot \nabla u =  \nu   \Delta u, \quad  u(\cdot,0) =  u_0.
\end{equation}
Under the above conditions, it is known that

\begin{displaymath}
\lim _{t \to \infty} \Vert u (t) \Vert _{L^2} = 0, \quad \lim _{t \to \infty} t^{n/4} \Vert u(t)\Vert _{L^{\infty}(\mathbb{R}^{n})} = 0,
\end{displaymath}
see Theorem 3.3 in Braz e Silva, Sch\"utz and P.R. Zingano \cite{MR3117150},  and also that
\begin{displaymath}
\| u(t) \| _{L^{\infty}(\mathbb{R}^{n})} \leq  K_{n} \| u_0 \|_{L^2}  \nu^{- n/4} t^{-n/4} \quad \forall  t > 0 
\end{displaymath}
for some constant  $K_{n}$ depending only on $ n$, see Theorem 3.2 in Braz e Silva, Sch\"utz e P.R. Zingano \cite{MR3117150} and Theorem 2.1 in Porzio \cite{MR2529737}. Conditions (\ref{eqn:evreg1}) and (\ref{eqn:evreg2}) hold for $t_{\ast} = 0$. From this point on, we restrict our attention to dimension $ n \geq 2 $.

We check now hypotheses (\ref{eqn:h2}), where $ \mathbb{G}(u) = \bb(u) \cdot \nabla u $. It clearly holds for $m = 0$, since 

\begin{displaymath}
\int_{\mathbb{R}^{n}}  u(x,t) \bb (u) \nabla u \, dx = \int_{\mathbb{R}^{n}} \nabla \cdot \mbox{\boldmath $D$}(u(x,t)) dx =0
\end{displaymath}
because $ u(t) \in L^{2}  \cap L^{\infty}(\mathbb{R}^{n}) $, for any $ t > 0 $, where

\begin{displaymath}
D(u) = \int_{0}^{u}  v \, \bb (v) \, dv.
\end{displaymath}
Checking (\ref{eqn:h2}) for $ m \geq 1 $ is more involved. We will need the following Lemma, whose proof can be found in page 273 in Moser \cite{MR199523} or page 70 in Klainerman \cite{MR544044}. 

\begin{Lemma}
\label{lemma-three-one}
Let  $m \geq  1 $,  $F \in C^{m}(\mathbb{R}) $ with $ F(0) = 0 $, and $u$ in $W^{m, p} (\mathbb{R}^{n})  \cap L^{\infty}(\mathbb{R}^{n}) $, for $ 1 \leq p \leq \infty$. Then $F(\mbox{u})$ is in $W^{m, p}(\mathbb{R}^{n})$ and 

\begin{displaymath}
\| D^{m} F(u) \|_{L^{p}(\mathbb{R}^{n})} \leq K (m,n,p) \, \| D^{m}u \|_{L^{p}(\mathbb{R}^{n})}
\end{displaymath}
where $K > 0$ depends only on $ m, n, p$ and
\begin{displaymath}
F_{\scriptstyle \ell} =  \sup \{ |F^{(\ell)}(v)| : |v| < \Vert u \Vert _{L^{\infty}(\mathbb{R}^{n})} \}
\end{displaymath}
for $ 1 \leq \ell \leq m $, where $ F^{(\ell)}$ denotes the derivative  of order $\ell $ of the function  $F$.
\end{Lemma}
Given $u \in \mathbb{R} $, it will be convenient to define
$\widetilde{\bb} (u), \widetilde{\BB} (u)$ in $\RR ^n$  through
\begin{displaymath}
\widetilde{\bb} (u) = \bb (u) - \bb (0), \qquad \widetilde{\BB} (u) = \int _0 ^u \left(\bb(v) - \bb (0) \right) \, dv, 
\end{displaymath}
which leads to

\begin{displaymath}
\widetilde{\mathbb{G}} (u) = \left(\bb (u) - \bb (0) \right) \cdot \nabla u = \mathbb{G} (u) - \bb (0) \cdot \nabla u = \nabla \cdot \widetilde{\BB} (u),
\end{displaymath}
where $ u = u(x,t)$ is the solution to (\ref{eqn:adv-diff}) under consideration. From now on, we will work with $\widetilde{\mathbb{G}} (u)$, as (\ref{eqn:h2}) is valid for $\mathbb{G} (u)$ if and only if is valid for $\widetilde{\mathbb{G}} (u)$, for any $m$.

The next lemma is proved by direct computation. 

\begin{Lemma}
\label{lemma-threetwo}
Let $v$ be a smooth  scalar function in $ \mathbb{R}^{n}$, and let $ D_{\ell} = \frac{\partial}{\partial x_{\ell}}$. Then 

\begin{displaymath}
D_{\ell} \widetilde{\BB} (v) = \widetilde{\bb} (v)  D_{\ell} v,
\end{displaymath}
and for $m \geq 2$

\begin{displaymath}
D_{\ell _1} D_{\ell _2} \cdots D_{\ell _m}  \widetilde{\BB} (v) = \widetilde{\bb} (v) D_{\ell _1} D_{\ell _2} \cdots D_{\ell _m} v + \sum _{j = 1} ^{m - 1} \left( \sum _{k=1} ^{a(m,j)} \mathbb{D} _{(k)} ^j \widetilde{\bb} (v) \cdot \mathbb{D} _{(k)} ^{m-j} v  \right),
\end{displaymath}
where 

\begin{displaymath}
a(m,j) = \frac{(m-1)!}{j! (m - 1 - j)!}, \, \mathbb{D} _{(k)} ^j = D_{ k_1} D_{k_2} \cdots D_{k_j}, \, \mathbb{D} _{(k)} ^{m - j} = D_{ k_{j + 1}} D_{k_{j + 2}} \cdots D_{k_m},
\end{displaymath}
with  $ \{\!\;\!\;\! k_{1}, k_{2}, \!\;\!...\!\;\!\;\!, k_{m} \} = \{\;\! \ell_{1}, \ell_{2}, \!\;\!...\!\;\!\;\!, \ell_{m} \!\;\!\;\!\}  $.
\end{Lemma}

\begin{Lemma}
\label{lemma-threethree}
Let $n \geq 2, m \geq 1$ and $F, G \in C^m (\RR)$, with $F(0) = G(0) = 0$. Let $\lambda_1, \lambda_2 \in \ZZ ^n$ non-negative multi-indices such that $|\lambda_1| + |\lambda_2| = m$. Then, for any $v \in H^{m+1} (\RR^n) \cap H^{n  - 2} (\RR^n) \cap L^{\infty} (\RR ^n)$, we have that $D^{\lambda_1} F(v) \cdot D^{\lambda_2} G(v)$ is in $L^2 (\RR ^n)$ and

\begin{displaymath}
\| D^{\lambda_1} F(v) \cdot  D^{\lambda_2}  G(v) \|_{L^2 (\RR ^n)} \leq C(m,n) \, \|\,v\, \|_{L^{2} }^{\frac{1}{2}} \|\,D^{n-2} v\, \|_{L^{2} }^{\frac{1}{2}}  \|\,D^{m+1} v\, \|_{L^{2} }
\end{displaymath}
where $ C(m,n) $ depends only on $m,n$ and

\begin{displaymath}
F_{\ell} = \sup\,\{|F^{(\ell)}(w)|: |w| < \Vert v \Vert _{L^{\infty} (\RR ^n)} \}, \, G_{\ell} = \sup\,\{|G^{(\ell)}(w)|: |w| < \Vert v \Vert _{L^{\infty} (\RR ^n)} \},
\end{displaymath}
for $1 \leq l \leq m$.
\end{Lemma}

\begin{Remark}
The assumption that $ v \in L^{\infty}(\mathbb{R}^{n}) $ in the previous Lemma  is needed when $ n = 4 $ and $ m = 1 $, but redundant in the other cases, since the inclusion $ H^{m+1}(\mathbb{R}^{n}) \cap   H^{n - 2}(\mathbb{R}^{n})  \subset L^{\infty}(\mathbb{R}^{n}) $ holds for $ m \geq 1 $ except in the single case  $ \:\! m = 1 $, $ n = 4 $. 
\end{Remark}

{\bf Proof:} The proof is analogous to Lemma 3.1 in Braz e Silva, J.P. Zingano and P.R. Zingano \cite{MR3907942}, to which we refer. $\Box$ 

\begin{Lemma}
\label{lemma-threefour}
Let $u(t) \in C([0, \infty), L^{2}(\RRn) )  \cap  L^{\infty}_{loc}( (0, \infty), L^{\infty}(\mathbb{R}^{n}))$
be a solution to (\ref{eqn:adv-diff}).  Then, for every $m \geq 1$, there exists a constant $K = K\!\;\!(m,n,\nu)$, depending only on $m, n, \nu $, $\bb$ and the size of  $\|u(t)\|_{L^{\infty}(\mathbb{R}^{n})} $ at, say,  $ t = 1$, such that
\begin{displaymath}
\| D^{m} u(t)\|_{L^{2} } \leq\, K (m,n,\nu) \, \| u_0 \|_{L^{2} } \quad  \forall \, t > m.
\end{displaymath}
\end{Lemma}

{\bf Proof:} From (\ref{eqn:adv-diff}) we have 

\begin{displaymath}
\| u(t) \|^{\:\!2} _{L^{2} } +\, 2 \nu \int_0 ^t  \| D u(\tau) \|^{2} _{L^{2} } \,d\tau = \|u_0 \|^{2} _{L^{2} }
\end{displaymath}
for all $t > 0 $. Then

\begin{displaymath}
\int_{0}^{1}  \| D u(\tau) \|^{2} _{L^{2} }  \,d\tau  \leq \frac{1}{2\nu} \| u_0 \|^{2} _{L^{2} } 
\end{displaymath}
so there exists  $ a_{1} \in (0, 1) $ such that 
\begin{equation}
\label{eqn:310b}
\| D u (a_{1}) \|^{2} _{L^{2} }  \leq  \frac{1}{2 \nu }  \| u_0 \|^{2} _{L^{2} }.
\end{equation}
Then, for all $t > a_1$ we have 
\begin{eqnarray*}
\| D  u (t) \|^{2} _{L^{2} }  & + &  \, 2 \nu \int_{a_1}^{t} \| D^{2} u (\tau) \|^{2} _{L^{2} } 
\,d\tau \leq \| D  u (a_{1}) \|^{ 2} _{L^{2} } \\  & + &  2 \int_{a_1}^{t} \| D^{2} u(\tau) \| _{L^{2} } \, \| D [\widetilde{\BB}(u)] \| _{L^{2} }  \,d\tau.
\end{eqnarray*}
Using Lemma \ref{lemma-three-one} we then obtain

\begin{eqnarray}
\label{eqn:310c}
\| D  u (t) \|^{2} _{L^{2} }   & + &  \nu \int_{a_1}^{t} \| D^{2} u (\tau) \|^{2} _{L^{2} } 
\,d\tau \leq \| D  u (a_{1}) \|^{ 2} _{L^{2} } \notag  \\ & + &   \frac{K_1}{\nu} \int_{a_1}^{t} \| D^{2} u(\tau) \| _{L^{2} }  \,d\tau.
\end{eqnarray}
Then, from (\ref{eqn:310b}) and (\ref{eqn:310c}) we obtain the result for $m = 1$ and also

\begin{displaymath}
\int_1 ^t \| D^{2} u(\tau) \|^{2} _{L^{2}} \,d\tau \,\leq  C(n,\nu) \, \| u_0 \|^{2} _{L^{2}}
\end{displaymath}
for all $t > 1 $ and some constant  $ C(n,\nu) > 0 $ that depends also on $\bb$ and  $ \| u(1)\|_ {L^{\infty}}$.
In particular,  we can find $a_{2} \in (1, 2) $ such that

\begin{displaymath}
\| D^{2} u(a_{2}) \|^{2} _{L^{2}} \leq C(n,\nu) \, \| u_0 \|^{2} _{L^{2}},
\end{displaymath}
starting with which we can consider the energy estimate

\begin{eqnarray*}
\| D^{2} u (t) \| ^{2} _{L^{2}} &+&\, 2\ \nu \int_{a_2}^{t} \| D^{3} u (\tau) \|^{2} _{L^{2}}  \,d\tau \,\leq  \| D^{2} u(a_{2}) \|^{2} _{L^{2}} \\&+& 2  \int_{a_2}^{\t} \|D^{3} u (\tau) \|_{L^{2}}  \, \| D^{2} [\widetilde{\BB} (u) ] \| _{L^{2}} \,d\tau.
\end{eqnarray*}
We now proceed as above to obtain the result for $m = 2$, and go on to prove the estimate for any $m$. $\Box$

We are ready now to show that solutions to (\ref{eqn:adv-diff}) satisfy (\ref{eqn:h2}) for any $ m \geq 0$. This is so,  because given $t \geq 1$, $ m \geq 1 $ and $\ell_{1}, \ell_{2}, \cdots , \ell_{m} \in \{1, 2, \cdots, n\} $, using Lemmas \ref{lemma-threetwo} and \ref{lemma-threethree} we obtain 

\begin{align}
\biggl|\ \int_{\mathbb{R}^{n}} D_{\ell_1} \cdots  D_{\ell_m} u(x,t) & \cdot  D_{\ell_1} \cdots D_{\ell_m} \mathbb{G}(u(x,t))  \,dx \biggr|  \leq    \|D^{m +1} u (t)\ \|_{L^{2}} \| D^{m} [\,\widetilde{\BB} (u(t))] \| _{L^{2}} \notag \\ & \leq  C(m,n) \, \| u (t) \|^{\frac{1}{2}} _{L^{2}} \| D^{n-2} u (t) \|^{\frac{1}{2}} _{L^{2}} \| D^{m +1} u (t) \|^{2} _{L^{2}} \notag
\end{align}
where $C(m,n)$ denotes some constant that depends on $m, n$ and 
\begin{displaymath}
b_{\ell} =  \max_{1 \leq \ell \leq m} \left\{|\bb ^{(\ell)} (v)|: |v| \leq \Vert u(1) \Vert _{L^{\infty}} \right\}.
\end{displaymath}
Then, Lemma \ref{lemma-threefour} leads us to the result.

Now, suppose 
\begin{displaymath}
\| u(t) \| _{L^2}  \leq  C_{0}  t^{-\,\alpha},  \qquad  \forall t > T_{0}
\end{displaymath}
for some constants $\alpha, C_0, T_{0} > 0$. Then, Theorem \ref{upper-bound-decay} leads us to 

\begin{displaymath}
\| D^{m} u(t) \| _{L^2}  \leq  C_{m} \nu^{-\frac{m}{2}} t^{- \alpha - \frac{m}{2}},  \qquad \,  \forall t > T_{m}
\end{displaymath}
for some  $ T_{m} > T_{0} $ that  can be chosen to depend only on $ \:\!m, \:\!\nu $, $ C_0, T_{0} $ and $\:\!\alpha $, and for some $C_m = C_m(\alpha, C_0) > 0$  \\

Finally, assume also the estimate

\begin{displaymath}
\| u(t) \|_{L^2} \geq  c_0  t^{- \eta},  \qquad  \forall  t >  t_{0}
\end{displaymath}
holds,  for some constants  $c_0, t_0, \eta > 0 $, with $ 0 < \alpha \leq \eta $. Then, from Theorem \ref{lower-bound-decay-first} and \ref{lower-bound-decay-second} we obtain, for some $c_m = c_m (m, \alpha, C_0, c_0)$ 

\begin{displaymath}
\| D^{m} u(t) \| _{L^2} \geq\: c_m \nu^{-\frac{m}{2}}  t^{-\eta - \frac{mq}{2}}, \qquad \forall t > t_{m}
\end{displaymath}
where $q = \frac{\eta}{\alpha}$ and $ t_{\!\:\!m} \!\:\!> t_{\!\:\!0} \!\;\!$ that can be chosen to depend only on $ \;\!m $, $ \nu $, $ t_0 $, $\!\:\! T_{0} $,  $ c_0 $, $ C_0, \:\!\alpha $ and $ \eta $.

\subsection{Incompressible MHD equations} Let $2 \leq n \leq 4$.  Given  $\uu_0, \bb_0 \in L^{2}_{\sigma}(\mathbb{R}^{n})$, let  $\uu, \bb  \in  C_{w}([0, \infty), L^{2}_{\sigma}(\mathbb{R}^{n})) \cap L^{2}((0,\infty), \dot{H}^{1}(\mathbb{R}^{n})) $ be any given Leray \cite{MR1555394} solution to the MHD equations 

\begin{displaymath}
\begin{dcases*}
\uu_t + \left(\uu \cdot \nabla \right) \uu + \nabla p = \mu \Delta \uu + \left( \bb \cdot \nabla \right) \bb, & $\nabla \cdot \uu = 0$ \\
\bb_t + \left(\uu \cdot \nabla \right) \bb  = \nu \Delta \bb + \left( \bb \cdot \nabla \right) \uu, & $\nabla \cdot \bb = 0$,
\end{dcases*}
\end{displaymath} 
with $\uu_0 = \uu ( 0)$, $\bb_0 = \bb ( 0)$. It is known that these solutions are eventually smooth, see Theorem 4.2 in Sermange and Temam \cite{MR716200} for $n = 2, 3$ and Melo, Perusato, Guterres and Nunes \cite{MR4128674}. If $n = 4$, the result holds provided $\Vert (\uu _0, \bb_0) \Vert _{\dot{H} ^1}$ is small enough, see Lemma 3.1 in  Braz e Silva, J.P. Zingano and P.R. Zingano \cite{MR3097244}. If we consider
\begin{equation}
\notag
\mathbb{G}(\uu, \bb) = \left( \begin{tabular}{c} $\mathbf{f}$ \\ $ \mathbf{g}$ \end{tabular} \right) = 
\left( \begin{tabular}{c}
$\left(\uu \cdot \nabla \right) \uu + \nabla p - \left( \bb \cdot \nabla \right) \bb$ \\
$\left(\uu \cdot \nabla \right) \bb  - \left( \bb \cdot \nabla \right) \uu$
\end{tabular}
\right),
\end{equation}
then (\ref{eqn:h2}) is valid for any $m \geq 1$, see Lemma 3.1 in Braz e Silva, J.P. Zingano and P.R. Zingano \cite{MR3907942}. Note that $\mathbf{f}, \mathbf{g}$  satisfy (\ref{eqn:h3}) for all $m \geq 1$  with $ \beta = 2 \alpha + \frac{n+2}{4}$.

We now assume that 

\begin{equation}
\label{eqn:five-three}
\| \left( \uu, \bb \right) (t) \|_{L^2}  \leq C_{0} \: t^{-\,\alpha},  \qquad \, \forall  t >  T_{0}
\end{equation}
for some  $ \alpha, C_0, T_{0} > 0 $. Then, using Theorem \ref{upper-bound-decay} we obtain that

\begin{displaymath}
\| D^{m} \uu (t) \| _{L^2} \leq C_m  \mu^{-\frac{m}{2}} t^{- \alpha - \frac{m}{2}},  \qquad \forall  t > T_{m}
\end{displaymath}
and
\begin{displaymath}
\| D^{m} \bb (t) \| _{L^2} \leq C_m  \nu^{-\frac{m}{2}} t^{- \alpha - \frac{m}{2}},  \qquad \forall  t > T_{m}
\end{displaymath}
where $C_m$ depends on $m$, $ \alpha $ and $C_0$ and for some  $ T_{m} > 0$ that depends on 
$m, \alpha, \mu, \nu$, $C_{0}$ and $T_{0} $.

If we assume \eqref{eqn:five-three} and the lower bound

\begin{equation}
\label{eqn:five-four}
\| \left( \uu, \bb \right) (t) \|_{L^2}  \geq c_{0} \: t^{-\,\alpha},  \qquad \, \forall  t >  t_{0},
\end{equation}
for some $c_0, t_0 >0$, using Theorem \ref{lower-bound-decay-first} we obtain

\begin{displaymath}
\| D^{m} (\uu,\bb) (t) \| _{L^2} \geq c_m  \,\gamma^{-\frac{m}{2}}\, t^{- \alpha - \frac{m}{2}},  \qquad \forall  t > t_{m}
\end{displaymath}
where $ \gamma =\min\{\mu,\nu\} $ and $c_m$ depends on $m$, $ \alpha $, $ c_0 $ and $C_0$, for some  $ t_{m} > 0$ that depends on 
$m, \alpha, \mu, \nu, c_0$, $C_{0}$, $t_0$ and $T_{0} $. 
\begin{Remark}
	As another illustration of Theorem \ref{lower-bound-decay-first}, let us assume \eqref{eqn:five-three} and that \eqref{eqn:five-four} applies
	only to $ \mbox{\boldmath $u$}(t) $, say,
	so that we have 
	\begin{equation*}
	\|  \uu  (t) \|_{L^2}  \geq c_{0} \: t^{-\,\alpha},  \qquad \, \forall  t >  t_{0},
	\end{equation*}
	for some $\!\;\!\;\! c_0, \!\;\!\;\!t_0 > 0 $, with no lower bound assumed for $ \mbox{\boldmath $b$}(t) $.  In this case, 
	Theorem \ref{lower-bound-decay-first} would give 
	\begin{equation}
	\notag\|\:\!D^{m} \mbox{\boldmath $u$}(t) \:\!
	\|_{L^{2}} \:\!\geq\,
	c_{m} \, \mu^{-\;\!\frac{m}{2}} \,
	t^{\:\!-\;\!\alpha\;\!-\;\!\frac{m}{2}}
	\!\;\!,
	\quad \;\;\,
	\forall \; t > t_{m} 
	\end{equation}
	for every $m$,
	with $ c_m \!\!\;\!\;\!= c_m(\alpha, c_0, C_0) $
	and $ t_m \!\!\;\!\;\!= t_m(\alpha, \mu, \nu, c_0, C_0, t_0, T_0) $. Or, assuming \eqref{eqn:five-three} and that, say,
	\begin{equation*}
	\|  b_i  (t) \|_{L^2}  \geq c_{0} \: t^{-\,\alpha},  \qquad \, \forall  t >  t_{0},
	\end{equation*}
	for some $1 \leq i \leq n$, then we would obtain from Theorem \ref{lower-bound-decay-first} that
	\begin{equation*}
	\| D^m  b_i  (t) \|_{L^2}  \geq c_{m}\,\nu^{-\frac{m}{2}} \: t^{-\,\alpha - \frac{m}{2}},  \qquad \, \forall  t >  t_{m},
	\end{equation*}
	for every $m$, 
	and so forth.
\end{Remark}

These lower bounds can be extended to the case where $\eta > \alpha$ by using Theorem \ref{lower-bound-decay-second}. For example, if 

\begin{equation}
\label{eqn:five-five}
\|  \uu  (t) \|_{L^2}  \geq c_{0} \: t^{-\,\eta}\,\,\,\,\mbox{and}\,\,\,\,\,\|  \bb (t) \|_{L^2}  \geq c_{0} \: t^{-\,\eta},   \qquad \, \forall  t >  t_{0}
\end{equation}
for some  $ \alpha, C_0, T_{0} > 0$, where $\eta = q \alpha$, with  $ 1 < q < \frac{3\alpha + \frac{n+2}{4}}{2 \alpha +1}$, then by Theorem \ref{lower-bound-decay-second} there exists $c_1 = c_1 (\alpha, c_0, C_0)$ such that 

\begin{displaymath}
\| D \uu (t) \| _{L^2} \geq c_1  \mu^{- \frac{1}{2}} t^{- \eta \,- \frac{q}{2}}  \qquad \, \forall t  >  t_{1}
\end{displaymath}
and
\begin{displaymath}
\| D \bb (t) \| _{L^2} \geq {c}_1  \nu^{- \frac{1}{2}} t^{- \eta \,- \frac{q}{2}}  \qquad \, \forall t  >  t_{1}, 
\end{displaymath}
for some  $ t_{1}> 0$ that depends on  $\alpha, q, \mu, \nu, c_0, C_{0}, t_0 $ and $T_{0} $.

\bibliographystyle{plain}
\bibliography{GuterresNichePerusatoZingano}{}

\end{document}